\documentclass[11pt]{article}
\usepackage{url}
\usepackage{amsmath} 
\usepackage{times}
\usepackage{amsfonts} 
\usepackage{graphics}
\usepackage{color}
\newtheorem{corollary}{Corollary}[section]

\def\ve{\varepsilon} 
\def\({\left(} 
\def\){\right)} 
\def\[{\left[} 
\def\]{\right]} 
\def\norm#1{\left\| #1 \right\|}

\newcommand{\closeproof}{\hfill $\Box$}
\def\norm#1{\left\| #1 \right\|}
\def\absval#1{\left| #1 \right|}

\newtheorem{theorem}{Theorem}[section]

\begin{document}

\title{Identification of a chemotactic sensitivity in a coupled system.} 
\author{K. Renee Fister\\ Department of Mathematics and Statistics \\
   Murray State University \\ Murray, KY 42071 \\ 
\texttt{renee.fister@murraystate.edu}
   \bigskip
   \\ Maeve L.  McCarthy \\ Department of Mathematics and Statistics \\
   Murray State University \\ Murray, KY 42071 \\ 
   \texttt{maeve.mccarthy@murraystate.edu}}
\date{\today}
\maketitle
\begin{abstract} 
Chemotaxis is the process by which cells behave in a way that follows the chemical gradient.  Applications to bacteria growth, tissue inflammation, and vascular tumors provide a focus on optimization strategies.  Experiments can characterize the form of possible chemotactic sensitivities.  This paper addresses the recovery of the chemotactic sensitivity from these experiments while allowing for nonlinear dependence of the parameter on the state variables. The existence of solutions to the forward problem is analyzed.  The identification of a chemotactic parameter is determined by inverse problem techniques. Tikhonov regularization is investigated and appropriate convergence results are obtained. Numerical results of concentration dependent chemotactic terms are explored. 
\end{abstract}



\section{Introduction} 
\label{sec: introduction}

Biological and ecological research has investigated cell migration.  To model cell migration, studies have been composed to include migration, diffusion, haptotaxis, and chemotaxis, 
\cite{Keller, Kell, Keller3, Oster, Anderson1, Anderson2, Dung}.  In this paper, the focus is chemotaxis.  Chemotaxis describes the movement of an organism and/or groups of cells that either move toward or away from a chemical or sensory stimulus.  In the early work by Keller and Segel \cite{Keller}, chemotactic responses of amoebae to bacteria is studied in a cellular slime mold.  Bacterial chemotaxis, which describes the ability of bacteria to move toward increased or decreased concentrations of attractants is analyzed at the macroscopic level through a microscopic model of individual cells, \cite{Erban, Seg}.  It was first observed by Engelmann \cite{Engelmann} in 1881.  For example, if {\em Salmonella typhimurium}, a strain of salmonella associated with meat and poultry products, is introduced to a petri dish filled with a nutrient, the bacteria will migrate outward, consuming the nutrient. As they consume the nutrient, they secrete a chemoattractant. After several days, the bacteria will have clustered in the areas of high chemical concentration. A structure of concentric rings is usually observed experimentally.  Chet and Mitchell's work \cite{Chet} describes patterns formed from {\it{E. coli}} movement toward amino acids.  Allweis et. al \cite{Allweiss} investigate {\it{Vibrio cholerae}} which are inhibited by a pepsin digest that reduces the possibility of the vibrios attaching to the intestinal wall.   Other authors \cite{Fisher, Alt} have addressed chemotaxis in immune cell motility which when combined with tumor morphology is hoped to provide new avenues of treatment strategies.  In addition, authors have analyzed chemotactic responses in ecology \cite{Lapidus} and have investigated mathematical issues for the existence of global solutions in multiple dimensions, \cite{Kowalczyk, Herrero, Hillen, Childress, Chalub, Aida, Jabbarzadeh, SenbaSuzuki, Horstmann, SenbaSuzuki2, GajewskiSkrypnik}.

Chemotaxis also arises in a variety of medical applications.  In particular, it has been studied in connection with myxobacteria \cite{Slius,Tosin}, leukocyte mobility in tissue inflammation \cite{Alt}, the migration of tumor cells towards bone \cite{Orr},  and other issues in morphogenesis \cite{Maini}. Another interesting problem involves the study of vascular tumors through angiogenesis.  Angiogenesis involves the formation of capillary networks of blood vessels that are vital for the growth of tumors.  Mathematical modeling of angiogenesis \cite{Anderson1, Anderson2, Bellomo, Corrias, Hillen2, Levine, Orr, Painter, Sleeman, Vel1, Vel2} has given new insight into tumor structure.  Normal tissue, lymphocytes, and other types of cells grow at the tumor site or are recruited through chemotaxis.  The need to identify the nature of this recruitment is at the heart of this paper.  The identification of a chemotactic term falls under the umbrella of an inverse problem.  In principle, we can measure certain characteristics of the tumor concentration and use mathematical techniques to recover the chemotactic term, in particular the chemotactic sensitivity, that is driving the tumor growth.  To our knowledge, this {\em inverse problem} approach has only been used in the analysis of chemotaxis models by Dolak-Stru{\ss} and K\"{u}gler \cite{DolakKugler} under the assumption that the chemical concentration is explicitly known. 


Since there are many applications in which chemotaxis arises, there are also different models of the chemotactic effect.
There have been many different expressions proposed that model chemotactic velocity, see Keller and Segel \cite{Keller3}, Lapidus and Schiller \cite{Lapidus2},  Ford and Lauffenburger \cite{Ford}, and Tyson et al. \cite{Tyson}.
 This velocity is used in a bacterial conservation equation in the formulation of a system of partial differential equations that governs the particular application. The chemotactic sensitivity determines the velocity.
Our goal is to develop a technique whereby the appropriate chemotactic sensitivity model, and hence chemotactic velocity, can be determined from available data. 

In particular, we consider a system of partial differential equations that was developed by Oster and Murray \cite{Oster} to model  the  pattern formation of cartilage condensation in a vertebrate limb bud. A similar system was studied by Myerscough et al., \cite{Myerscough}. The numerical solution of similar systems were recently studied by Tyson et al. \cite{Tyson2} and by Nakaguchi and Yagi \cite{Nakaguchi}. Work by Fister and McCarthy \cite{mcfist} has shown that the system of partial differential equations can in fact be controlled theoretically through the introduction of a mechanism controlling the number of cells being generated. Simulations provide optimal drug treatment programs for patients to facilitate the rebuilding of cartilage or the reduction of cancerous tumors. The chemotactic sensitivity in \cite{mcfist} was known and the control parameter was a harvesting term. Our goal in this work is to {\em identify} the chemotactic sensitivity.

The paper is organized into six sections.  In section two, the existence of the forward problem is proven.  In section three, identifiability of the chemotactic sensitivity is established using the weak formulation of the state problem.  In section four, Tikhonov regularization is used to approximate the solution through the use of minimization arguments.  The rate of convergence of the approximate minimizer of the chemotactic sensitivity to the true parameter follows next.  In section five, numerical experiments provide graphical depictions of the accuracy of the recovery of the parameter.  In section six, conclusion remarks are made.

\section{Forward Problem}
In this model, $u(x,t)$ and $c(x,t)$ represent the concentration of the cells and the chemoattractant, respectively. The cells 
and the chemoattractant are governed by a convection-diffusion equation and a reaction-diffusion equation as
\begin{eqnarray} 
&u_t = M\Delta u - \nabla \cdot (\chi(u,c) u  \nabla c) & \mbox{in $\Omega \times (0,T)$} 
\\ \nonumber 
&c_t= D\Delta c + {bu\over u+h} -  \mu c &\\ \nonumber 
&u(x,0) =u_0(x), \quad c(x,0) = c_0(x) & \mbox{for $x\in \Omega$}\\ 
&\frac{\partial u}{\partial \nu}= \frac{\partial c}{\partial \nu}= 0 & \mbox{on $\partial\Omega\times (0,T)$}
\nonumber 
\end{eqnarray}
where $\nu$ is the outward unit normal. $M$ and $D$ represent the diffusion 
coefficients of the cells and the chemoattractant.  
 The Michaelis-Menten term, ${bu\over u+h}$, represents a 
response of the chemoattractant to a maximum carrying capacity or saturation 
rate, assuming $b,h>0$. We incorporate a decay term where $\mu$ denotes the 
degradation rate.
We assume that there is no flux of the 
concentrations across the boundary, and that the initial concentrations for the cells and chemoattractant are $u_0(x)$ and $c_0(x),$ respectively.

Here, $\chi(u,c)$ is the chemotactic sensitivity which monitors the chemical gradient attraction 
of the cells.  It is this term that we seek to identify. 
In \cite{Oster,Myerscough,mcfist} the term $\chi(u,c)$ is simply  a constant. %
More generally, $\chi(u,c)$ is a linear function of $u$  in \cite{VelazquezA,VelazquezB,HillenPainter, PainterHillen}, while in \cite{Keller,Keller3,Lapidus2,Ford,Tyson} it is a nonlinear function of $c.$
We assume henceforth that the chemotactic sensitivity has the form  $\chi(u,c)=a(c)$ and  is a bounded function. We restrict our analysis to the dimensionless system
\begin{eqnarray} 
\label{forward} 
&u_t = M\Delta u - \nabla \cdot (a(c)u\nabla c) & \mbox{in $\Omega \times (0,T)$} 
\\ \nonumber 
&c_t= D\Delta c + {u\over u+1} -  c &\\ \nonumber 
&u(x,0) =u_0(x), \quad c(x,0) = c_0(x) & \mbox{for $x\in \Omega$}\\ 
&\frac{\partial u}{\partial \nu}= \frac{\partial c}{\partial \nu}= 0 & \mbox{on $\partial\Omega\times (0,T).$}
\nonumber 
\end{eqnarray}
We will establish a technique for the identification of $a(c) \in \mathcal{A}$ where
\begin{eqnarray*}
\mathcal{A} &=& \left\{ a \in H^1(I) : \, 
\norm{\frac{\partial a}{\partial c} (c_1) - \frac{\partial a}{\partial c} (c_2) }_{L^2(I)}\leq K\norm{c_1-c_2}_{L^2(I)}
\right\}
\label{A-defn}
\end{eqnarray*}
Observe that, with the available data, we can only expect to recover $a(c)$ on the interval $I=\left[ c_{\min},c_{\max}\right].$ 
The Lipschitz condition on the derivative of the chemotactic parameter is quite reasonable, since the chemotactic parameter has a rate of change that is bounded for bacteria growth, \cite{Ford}. 

In order to prove identifiability and to establish the rate of convergence ot our method, we will need to establish existence of a solution of (\ref{forward}). Using the standard notation $H^k(\Omega)$ to represent the Sobelev space $W^{k,2}(\Omega),$ 
let $H^{k+\theta}(\Omega)$ denote the intermediate space between $H^k(\Omega)$ and 
$H^{k+1}(\Omega)$ for any $0<\theta<1.$ Let $D$ be an interval in $[0,\infty).$ The space 
$L^p(D;X)$ is the $L^p$ space of measurable functions in $D$ with values in the Banach space 
$X.$ The space $C^m(D;X), m=0,1,2,\ldots$ is the space of $m-$times continuously 
differentiable functions in $D$ with values in $X,$ while the space 
$C^{\theta}(D;X),0<\theta<1$ is the space of H\"{o}lder-continuous functions in $D$ with 
values in $X.$ 
\begin{theorem} 
\label{existencet} 
If $u_0, c_0 \in H^{1+\ve} (\Omega)$ for $0<\ve\le 1$, 
and $u_0(x)\ge 0$, $c_0(x) \ge \overline c_0 > 0$ on $\overline \Omega$, 
then a real unique local solution $u, c$ of (\ref{forward}) exists on an 
interval 
$[0,T]$ such that
$$ 
u, c \in C^{\eta} \([0,\infty) ; H^{1+\ve_1} (\Omega) \) \cap C \([0,T) ; H^2 
(\Omega) \) \cap C^1 \([0,T) ; L^2 (\Omega) \)$$
with $0<\ve_1<\min(\ve,{1\over2})$ and $0<\eta<\min \( 
{\ve-\ve_1\over2},{1-2\ve_1\over4} \)$. The solution satisfies the lower 
bounds 
$$ 
u(x,t) \ge 0, \quad c(x,t) \ge \overline c_0 e^{-t} \quad \mbox{ on [0,T].} 
$$ 
\end{theorem}

{\bf Proof:} 
Let $X=L^2(\Omega) \times L^2(\Omega)$ and $Z= H^{1+\ve}(\Omega) \times H^{1+\ve}(\Omega)$. 
The system (\ref{forward}) can be formulated as an abstract quasilinear 
equation 
\begin{eqnarray*} 
{dv\over dt} &+& A(v)v=f(v) \qquad \mbox{$0<t<\infty$} \\ 
v(0) &=& v_0 
\end{eqnarray*}
on the Banach space $X.$

Let 
$$v = \left( \begin{array}{r} u \\ c\end{array}\right), \quad \hat v = \left( \begin{array}{r} \hat u \\ \hat c\end{array}\right), 
\quad v_0 = \left( \begin{array}{r} u_0 \\ c_0\end{array}\right).$$
Clearly $v_0\in Z$. 

We define $A(v)$ to be the linear operator in $X$ such  that 
$$ 
A(v) \hat v =\nabla \cdot \left( \begin{array}{cc} - M  &    a(c)u   \\  
0 & - D    \end{array} \right)
 \left( \begin{array}{c}  \nabla\hat u \\  \nabla\hat  c  \end{array} \right) 
+
 \left( \begin{array}{cc}  M  &  0   \\  
0 &  1  \end{array} \right) 
 \left( \begin{array}{c} \hat u \\ \hat c  \end{array} \right) 
$$ 
with domain 
$$ 
D\(A(v)\) = \biggl\{ \hat v \in H^2(\Omega)\times H^2(\Omega) ; 
 \frac{\partial \hat u}{\partial n } =  \frac{\partial \hat c}{\partial n } = 0
 \mbox{ on }  \partial\Omega \biggr\}. 
$$ 
Let the vector $f(v)$ be the function 
$$ 
f(v) =  \left( \begin{array}{c} Mu \\  \frac{u}{u+1}  \end{array} \right) . 
$$ 
Since $f(v)$ is Lipschitz, application of  Yagi's work \cite[Thm 2.1 and 3.4]{Yagi} 
 yields our result.   (See Appendix A for statements of Yagi's results.)
\closeproof

\section{Inverse Problem Statement and Identifiability}
\label{subsec:Identify}
In this section, we begin by establishing the identifiability of the  parameter $a(c)$ from the available
data $u(x,t)$ and $c(x,t)$ almost everywhere in $W=L^2((0,T), H^1(\Omega))$.
Note that, in order for chemotaxis to be observed biologically, cells must be present and a chemical gradient must exist. This means that $u(x,t)$ and $\nabla c(x,t)$ must be nonzero for a measurable subset of $\Omega\times(0,T).$

We denote by $(u_a, c_a)$ and $(u_b,c_b)$ the solution pairs of (\ref{forward}) with chemotactic sensitivities $a(c)$ and $b(c),$ respectively. 

\begin{theorem}
\label{identification}
Let $(u_a, c_a)$ and $(u_b, c_b)$ both be solutions in  $W\times  W$
of the direct problem (\ref{forward}) corresponding to $a(c_a)$ and $b(c_b)$.  If $u_a=u_b$ and $c_a=c_b$ 
almost everywhere in $\Omega\times[0,T]$, then $a(c)=b(c)$.
\end{theorem}
{\bf{Proof:}} We consider the weak form of the first equation of the direct problem (\ref{forward}) for $(u_a,c_a)$ and $(u_b,c_b)$ and subtract them.
\begin{eqnarray*}
\label{weakdiff}
\int_0^T  \int_{\Omega}  \frac{\partial}{\partial t}(u_a-u_b) \phi \ dt \,dx +M \int_0^T \int_{\Omega}(\nabla(u_a-u_b))\nabla \phi \ dx\, dt \\ 
= - \int_0^T \int_{\Omega} [a(c_a)u_a\nabla {c_a} - b(c_b) u_b\nabla {c_b}] \nabla{\phi} \ dx\, dt \\
\end{eqnarray*}
Since $u_a(x,t)=u_b(x,t)$ and $c_a(x,t)=c_b(x,t)$ a.e., 
this reduces to 
$$  \int_0^T \int_{\Omega} [a(c_a) - b(c_a) ] u_a\nabla {c_a} \nabla{\phi} \ dx\, dt  =0.$$
By definition, $\phi$ is in $W=L^2((0,T);H^1({\Omega}))$, then we can choose $\phi(x,t)=c_a(x,t)$.  Hence,
$$  \int_0^T \int_{\Omega} [a(c_a) - b(c_a) ] u_a \left( \nabla {c_a} \right)^2 \ dx\, dt  =0.$$
Our existence result,  Theorem \ref{existencet},  says that $u_a \geq 0.$ 
We also employ our biological assumptions that $u_a \neq 0$ and $\nabla c_a \neq 0$ on a measurable subset of  $\Omega\times(0,T).$
Thus
$$a(c) = b(c)$$
almost everywhere.
\closeproof

\section{Output Least Squares and Tikhonov Regularization}
We wish to identify a function $a(c)\in\mathcal{A}$ from noisy measurements $(z_u,z_c)$ of $(u_a,c_a).$
Recall that 
\begin{eqnarray*}
\mathcal{A} &=& \left\{ a \in H^1(I) : \, 
\norm{\frac{\partial a}{\partial c} (c_1) - \frac{\partial a}{\partial c} (c_2) }_{L^2(I)}\leq K\norm{c_1-c_2}_{L^2(I)}
\right\}
\end{eqnarray*}
where $I=\left[ c_{\min},c_{\max}\right].$ 

We define 
\begin{equation}
F(a) \equiv ( u_a(x,t), c_a(x,t) )
\label{Fdefn}
\end{equation}
with
$$F : \mathcal{A} \rightarrow W\times W.$$
In the presence of perfect data $(z_u,z_c)$, we would solve the non-linear ill-posed problem
\begin{equation}
F(a^0)=(z_u,z_c)  \label{nonlinear-eqn}
\end{equation}
where $(u_{a^0},c_{a^0})$ is the solution of the direct problem with $a=a^0.$
To do this using Tikhonov regularization would involve approximating the solution by minimizing
\begin{equation}
 \label{minimization_problem}
\min_{a \in \mathcal{A}}  \norm{F(a)-(z_u,z_c)}^2_{W\times W} +\alpha \norm{a-a^*}^2_{L^2(I)}
\end{equation}
where $\alpha>0$ is a small parameter and $a^*$ is an {\em a priori} guess of the true solution $a^0.$
In real applications, measurement errors mean that exact data is not available. Noisy data is assumed to have an error level $\delta,$ which means that
\begin{equation}
\int_{0}^T \norm{u-z_u^\delta}^2_{L^2(\Omega)}  \ dt \leq \delta^2, \qquad
\int_{0}^T \norm{c-z_c^\delta}^2_{L^2(\Omega)}  \ dt \leq \delta^2
\label{noise-defn}
\end{equation}

We assume attainability of a true solution, i.e. if $(z_u,z_c)\in W\times W$ there exists $a^0 \in \mathcal{A}$ such that
\begin{equation}
F(a^0)=(z_u,z_c).
\end{equation}
In the presence of noisy data $(z_u^\delta,z_c^\delta), $ the minimizer $a_\alpha^\delta \in \mathcal{A}$ of (\ref{minimization_problem}) 
minimizes
\begin{eqnarray}
J_\alpha(a) & \equiv &  \norm{F(a)-(z_u^\delta,z_c^\delta)}^2_{W\times W} +\alpha \norm{a-a^*}^2_{L^2(I)} \nonumber \\
&=&\int_{0}^T \norm{u-z_u^\delta}^2_{L^2(\Omega)}  dt  +
\int_{0}^T \norm{c-z_c^\delta}^2_{L^2(\Omega)}  dt  \nonumber \\
&& +\alpha \norm{a-a^*}^2_{L^2(I)} 
\label{tikobj} 
\end{eqnarray}
for appropriate choices of $a\in \mathcal{A}$ and $\alpha.$

We begin by establishing the weak-closedness of the map $F(a).$

\begin{theorem}
If $a_n \rightharpoonup a_* \in \mathcal{A}$ then $u_{a_n} \rightharpoonup u_{a_*}$ and $c_{a_n} \rightharpoonup c_{a_*}$ 
in $W.$
\end{theorem}
{\bf Proof:}  Here, we give the outline of the proof and refer the reader to \cite{Soil} for details.  
Using that the solution to the state system (\ref{forward}) is unique, one can define $u_{a_n}=u(a_n)$ and $c_{a_n}=c(a_n)$.  
A transformation involving $e^{-\lambda t}$ times each component of the solution pair is made with $\lambda$ to be chosen 
in order to obtain the boundedness of the solution in $W$.   The weak definition of the solution associated with the 
transformed $u_{a_n}$ and $c_{a_n}$ in equation (\ref{weakdiff}) is analyzed via Cauchy's inequality and the boundedness of the coefficients.  
Using the boundedness (independent of $n$) of the solution pairs, subsequences are extracted that converge weakly to $u_*$ and $c_*$.  
Lastly,  comparison results are used  so that one can pass to the limit in the weak formulation of the solution to show that 
$u_*=u_{a*}$ and $c_*=c_{a*}$.
\qquad
\closeproof

Existence of a minimizer $a_\alpha^\delta$ now follows from the lower semi-continuity of  $J_\alpha(a).$
\begin{corollary}
For any data $(z_u^\delta,z_c^\delta) \in W\times W,$ a minimizer $a_\alpha^\delta$ of (\ref{tikobj}) exists.
\end{corollary}
Continuous dependence on the data $(z_u^\delta,z_c^\delta)$ for fixed $\alpha,$ and the convergence of $a_\alpha^\delta$
toward the true parameter $a^0$ as the noise level $\delta$ and the regularization parameter $\alpha$ go to zero also follow from standard results \cite{Seidman}.
\begin{corollary}
For fixed $\alpha,$ the minimizers depend continuously on the data $(z^\delta_u,z^\delta_c).$
If $\alpha(\delta)$ satisfies
$$\alpha(\delta)\rightarrow 0, \qquad \delta^2/\alpha(\delta)\rightarrow 0 \qquad \mbox{ as }\delta \rightarrow 0$$
then $$\lim_{\delta\rightarrow 0} \norm{a_{\alpha}^{\delta}-a^0}_{L^2(I)}=0.$$
\end{corollary}

\subsection{Convergence Rate Analysis}

Although we have noted (without proof) the convergence of the minimizer $a_\alpha^\delta$ to the true parameter $a^0,$ the rate of convergence
may be arbitrarily slow. We wish to determine a source condition
that will guarantee a certain  rate of convergence.
Even when our regularization parameter $\alpha$ is comparable to our noise level $\delta,$ 
such a source condition will require assumptions involving $u$ and $a^0-a^*.$

Recall that we seek to solve the nonlinear problem (\ref{nonlinear-eqn}), $F(a) =(z_u,z_c),$ where $F(a)\equiv (u_a,c_a).$ 
The true solution is $a^0,$ and $a^*$ is an {\em a priori} guess. 
In order to apply the theory of Engl, Hanke, Kunisch and Neubauer \cite{Engl1,Engl} , we must establish the following:
\begin{itemize}
\item $F$ is Frechet differentiable,
\item $F^\prime$ is Lipschitz with $\norm{F^\prime(a)-F^\prime(b)} \leq \gamma \norm{a-b},$
\item there exists $w$ satisfying the source condition $ a^\dagger - a^* = F^\prime(a^\dagger)^*w,$
\item $\gamma \norm{w}< 1.$ 
\end{itemize}
In practice, although computing $F^\prime$ and $(F^\prime)^*$ is not difficult, it can be quite tricky to establish the Lipschitz condition on $F^\prime$ with our system of coupled nonlinear partial differential equations.
Instead, our approach involves developing a {\em source condition} without imposing differentiability constraints on $F.$ Thus we establish 
$O(\sqrt{\delta})$ convergence. This technique is also found in the work of Engl and K\"{u}gler, \cite{EnglKugler}.

\begin{theorem}
Suppose that there exists a function $w\in L^2((0,T);H^1(\Omega))$ satisfying
$$w(x,0)=w(x,T)=0, \quad \Delta w \in   L^2((0,T);L^2(\Omega))$$
such that for any $\Psi \in \mathcal{A}$
$$  \left<a^0-a^*,\Psi \right>_{L^2(I)} = 
\int_{0}^T \int_{\Omega} \Psi(c_{a^0}) u_{a^0} \nabla c_{a^0} \cdot \nabla w \ dx \ dt.$$
If $\alpha \sim \delta$ then
\begin{eqnarray*}
\int_{0}^T \norm{u_{a_\alpha^\delta}-z_u^\delta}^2_{L^2(\Omega)} + \norm{c_{a_\alpha^\delta}-z_c^\delta}^2_{L^2(\Omega)}  \ dt  &=& O(\delta^2) \\
\end{eqnarray*}
and
\begin{eqnarray*}
\norm{a^0- a_\alpha^\delta}_{L^2(I)}  &=&  O(\sqrt{\delta}).
\end{eqnarray*}
\label{convergence-theorem}
\end{theorem}
{\bf Proof:}
For clarity, we briefly describe the techniques used in this proof. Using that a minimizer to $J_{\alpha}(a)$ exists, we obtain an upper bound in terms of the error level $\delta$ and the norm of the difference in the minimizer and optimal $a$ values. We then use our source condition with the weak formulation of the cell and chemical differential equations to obtain a representation of the inner product of the appropriate differences of the approximating minimizers. This allows us to bound $J_{\alpha}(a)$. Specifically, we use Triangle and Young's inequalities to bound the time and spatial derivatives of the differences in the state variables. Integration by parts and H\"{o}lder's inequality enable us to successfully bound the spatial derivatives of the states in terms of the states themselves. Using the assumptions from $\cal A$ and choosing $\epsilon$ sufficiently small, we can obtain the error of order $\sqrt{\delta}$ with $\alpha \sim \delta$.

Since $a_\alpha^\delta$ is a minimizer of $J_\alpha(a),$ we have $J_\alpha(a_\alpha^\delta) \leq J_\alpha(a^0).$ Using our definition of noise level 
(\ref{noise-defn}), we find that
\begin{eqnarray}
\int_{0}^T \norm{u_{a_\alpha^\delta}-z_u^\delta}^2_{L^2(\Omega)}  + 
\norm{c_{a_\alpha^\delta}-z_c^\delta}^2_{L^2(\Omega)}   \ dt  +\alpha \norm{a_\alpha^\delta-a^*}^2_{L^2(I)} 
\nonumber\\
 \leq 
\int_{0}^T \norm{u_{a^0}-z_u^\delta}^2_{L^2(\Omega)}  + \norm{c_{a^0}-z_c^\delta}^2_{L^2(\Omega)}   \ dt 
 +\alpha \norm{a^0-a^*}^2_{L^2(I)}
\nonumber\\
 \leq 
 2\delta^2 +\alpha \norm{a^0-a^*}^2_{L^2(I)}.
\label{source-ineq1}
\end{eqnarray}
Adding $\alpha \norm{a^0-a_\alpha^\delta}^2_{L^2(I)}$ to both sides of the inequality and using inner product properties yields
\begin{eqnarray}
\int_{0}^T \norm{u_{a_\alpha^\delta}-z_u^\delta}^2 + \norm{c_{a_\alpha^\delta}-z_c^\delta}^2  \ dt 
 +\alpha \norm{a^0- a_\alpha^\delta}^2_{L^2(I)} \nonumber \\
 \leq  2\delta^2 +2 \alpha \left<a^0-a^*,a^0-a_\alpha^\delta\right>_{L^2(I)}.
\label{source-ineq2}
\end{eqnarray}
Observe that our source condition
$$ \left<a^0-a^*,\Psi \right>_{L^2(I)} = 
\int_{0}^T \int_{\Omega} \Psi(c_{a^0}) u_{a^0} \nabla c_{a^0} \cdot \nabla w \ dx \ dt $$
with $\Psi=a^0-a_\alpha^\delta,$ together with the weak forms of the cell equation in the forward problem (\ref{forward}) for $a^0$ and $a_\alpha^\delta$ is
\begin{eqnarray*}
\left<a^0-a^*,a^0-a_\alpha^\delta\right>_{L^2(I)} =
 \int_{0}^T \int_{\Omega} \left( u_{a^0} -u_{a_\alpha^\delta}\right)_t w \ dx \ dt \\
+ \int_{0}^T \int_{\Omega} \nabla \left( u_{a^0} -u_{a_\alpha^\delta}\right) \cdot \nabla w\ dx \ dt \\
 +\int_{0}^T \int_{\Omega} \left[ a_\alpha^\delta(c_{a_\alpha^\delta})u_{a_\alpha^\delta} \nabla c_{a_\alpha^\delta} 
 -a_\alpha^\delta(c_{a^0}) u_{a^0} \nabla c_{a^0} \right] \cdot \nabla w\ dx \ dt
\end{eqnarray*}
and (\ref{source-ineq2}) becomes
\begin{eqnarray}
\int_{0}^T \norm{u_{a_\alpha^\delta}-z_u^\delta}^2_{L^2(\Omega)}  + \norm{c_{a_\alpha^\delta}-z_c^\delta}^2_{L^2(\Omega)}   \ dt 
 +\alpha \norm{a^0- a_\alpha^\delta}^2_{L^2(I)} \nonumber \\
 \leq  2\delta^2 +
 2\alpha\int_{0}^T \int_{\Omega} \left( u_{a^0} -u_{a_\alpha^\delta}\right)_t w \ dx \ dt \nonumber \\
+ 2\alpha M \int_{0}^T \int_{\Omega} \nabla \left( u_{a^0} -u_{a_\alpha^\delta}\right) \cdot \nabla w\ dx \ dt \nonumber \\
 +2\alpha\int_{0}^T \int_{\Omega} \left[ a_\alpha^\delta(c_{a_\alpha^\delta})u_{a_\alpha^\delta} \nabla c_{a_\alpha^\delta} 
 -a_\alpha^\delta(c_{a^0}) u_{a^0} \nabla c_{a^0} \right] \cdot \nabla w\ dx \ dt.
\label{source-ineq3}
\end{eqnarray}
We bound each integral in (\ref{source-ineq3}) separately using Triangle and Young's inequalities. For the estimates of $I_1,I_2$ and $I_3,$ we refer the reader to the Appendix for some of the details of the bounds used. We find that
\begin{eqnarray*}
\absval{I_1} &=&
 \absval{\alpha\int_{0}^T \int_{\Omega} \left( u_{a^0} -u_{a_\alpha^\delta}\right)_t w \ dx \ dt} \\
&\leq & 
\varepsilon \int_{0}^T \norm{u_{a_\alpha^\delta}-z_u^\delta}^2_{L^2(\Omega)} \ dt 
+ \frac{\alpha^2}{2\varepsilon}   \int_{0}^T \norm{w_t}^2_{L^2(\Omega)} \ dt 
+ \varepsilon \delta^2
\end{eqnarray*}
and
\begin{eqnarray*}
\absval{I_2} &=&
\absval{\alpha M \int_{0}^T \int_{\Omega} \nabla \left( u_{a^0} -u_{a_\alpha^\delta}\right) \cdot \nabla w\ dx \ dt} \\
& \leq & \varepsilon M^2\delta^2 
+ \frac{\alpha^2}{2\varepsilon}   \int_{0}^T \norm{\Delta w}^2_{L^2(\Omega)} \ dt 
+ \varepsilon M^2 \int_{0}^T \norm{u_{a_\alpha^\delta}-z_u^\delta}^2_{L^2(\Omega)} \ dt 
\end{eqnarray*}
where $\varepsilon$ is an arbitrary  parameter resulting from the use of Young's inequality. 
We utilize the assumptions 
$$\eta_1 \leq \frac{\partial a}{\partial c} \leq \hat{\eta}_1, \qquad 
\norm{\frac{\partial a}{\partial c} (c_1) - \frac{\partial a}{\partial c} (c_2) }_{L^2(\Omega)}\leq K\norm{c_1-c_2}_{L^2(\Omega)},$$
Green's Theorem, the boundary conditions, and H\"{o}lder's inequality to obtain the estimate,
%
%
%
\begin{eqnarray*}
\absval{I_3} &=&
\absval{\alpha\int_{0}^T \int_{\Omega} \left[ a_\alpha^\delta(c_{a_\alpha^\delta})u_{a_\alpha^\delta} \nabla c_{a_\alpha^\delta} 
 -a_\alpha^\delta(c_{a^0}) u_{a^0} \nabla c_{a^0} \right] \cdot \nabla w\ dx \ dt} \\
&\leq &
\varepsilon \left[ \hat{\eta}_1^2  \int_{0}^T \norm{u_{a_\alpha^\delta}-z_u^\delta}^2_{L^2(\Omega)} \ dt  + 
\norm{u}_{L^\infty(\Omega)}^2 K^2  \int_{0}^T \norm{c_{a_\alpha^\delta}-z_c^\delta}^2_{L^2(\Omega)} \ dt \right] \\
&&
+ \frac{\alpha^2}{\varepsilon} \int_{0}^T \norm{\Delta w}^2_{L^2(\Omega)} \ dt  
+ \varepsilon \hat{\eta}_1^2\delta^2
+ \varepsilon \norm{u}_{L^\infty(\Omega)}^2K^2\delta^2 \\
&&+ \varepsilon \norm{\nabla u}_{L^\infty(\Omega)}^2 K^2 \int_{0}^T \norm{c_{a_\alpha^\delta}-z_c^\delta}^2_{L^2(\Omega)} \ dt  + 
\varepsilon \hat{\eta}_1^2 \mu^2  \int_{0}^T \norm{u_{a_\alpha^\delta}-z_u^\delta}^2_{L^2(\Omega)} \ dt \\
&& + \frac{\alpha^2}{\varepsilon} \int_{0}^T \norm{\nabla w}^2_{L^2} \ dt  
+ \varepsilon  \norm{\nabla u}_{L^\infty(\Omega)}^2 K^2 \delta^2
+ \varepsilon \hat{\eta}_1^2\mu^2\delta^2.
\end{eqnarray*}

Grouping terms and relabeling constants, we have
\begin{eqnarray*}
&&
\int_{0}^T \norm{u_{a_\alpha^\delta}-z_u^\delta}^2_{L^2(\Omega)} + \norm{c_{a_\alpha^\delta}-z_c^\delta}^2_{L^2(\Omega)} \ dt 
 +\alpha \norm{a^0- a_\alpha^\delta}^2_{L^2(I)}\\
&\leq &  2\delta^2 +2 C_1 \varepsilon \delta^2 + 
2\varepsilon C_2 \int_{0}^T \norm{u_{a_\alpha^\delta}-z_u^\delta}^2_{L^2(\Omega)} \ dt \\
&+&
2\varepsilon C_3 \int_{0}^T \norm{c_{a_\alpha^\delta}-z_c^\delta}^2_{L^2(\Omega)} \ dt \\
&+&
\frac{\alpha}{\varepsilon} \int_{0}^T \left( \norm{w_t}^2_{L^2(\Omega)} +3\norm{\Delta w}^2_{L^2(\Omega)} + 2\norm{\nabla w}^2_{L^2(\Omega)} \right) \ dt. 
\end{eqnarray*}
If $\varepsilon$ is chosen to be sufficiently small, 
then for the choices $\alpha \sim \delta$  we obtain
$$\int_{0}^T \norm{u_{a_\alpha^\delta}-z_u^\delta}^2_{L^2(\Omega)} + \norm{c_{a_\alpha^\delta}-z_c^\delta}^2_{L^2(\Omega)}  \ dt =O(\delta^2)$$
and
$$\norm{a^0- a_\alpha^\delta}^2_{L^2(I)}  =O(\sqrt{\delta}). $$
\closeproof


\section{Numerical results}
\label{numerical-section}

In order to demonstrate the effectiveness of Tikhonov regularization for this application, we consider several examples. 

All computations were carried out in MATLAB. The Tikhonov functional
\begin{equation}
J_\alpha(a)
= \int_{0}^T \int_{\Omega} \left( |u-z_u^\delta|^2 + |c-z_c^\delta|^2  \right) \ dx \ dt 
 +\alpha \norm{ a-a^*}^2_{L^2(I)} 
\end{equation}
was minimized using {\tt lsqnonlin}, a MATLAB implementation of the Levenberg-Marquardt method with line search \cite{Levenberg,Marquardt}. Although it was not tractable to do so in the convergence analysis, a gradient based algorithm is appropriate here because computing the gradient and its adjoint is straightforward.

We restrict our discussion to $\Omega=[0,1].$ Recall that $z_u^\delta$ and $z_c^\delta$ represent noisy data and $a^*$ represents an {\em a priori} guess of the chemotactic sensitivity $a.$
Cell and chemoattractant concentration data on $\Omega=[0,1]$ was generated using {\tt pdepe} with high accuracy.
During the computation of $J_\alpha(a),$ cell and chemoattractant concentrations $u(x,t)$ and $c(x,t)$ associated with a particular $a$ were computed
using {\tt pdepe} with moderate accuracy over coarser space and  time meshes than those used to simulate data. 

Since {\tt lsqnonlin} requires objective functions of the form 
$$\frac{1}{2} \norm{F}_2^2 =\frac{1}{2}\sum_{k} f_k^2(x),$$
we approximated the first two terms of $J_\alpha(a)$ by
\begin{eqnarray*}
\sum_{j=1}^M \left\{ 
\sum_{i=1}^N \left[u(x_i,t_j)-z_u(x_i,t_j)\right]^2
+\sum_{i=1}^N \left[c(x_i,t_j)-z_c(x_i,t_j)\right]^2 
 \right\}  \left( \Delta x\right) \left(\Delta t\right) \\
\end{eqnarray*}
where $x_i=i(\Delta x)$ for $i=0,\ldots, N$ with $\Delta x = 1/N,$ and $t_j= j(\Delta t)$ for $j=0,\ldots, M$ with $\Delta t = \sigma/M.$

%

We approximate $a(c)$ by
\begin{equation}
a(c) =\sum_{k=1}^L a_k\phi_k(c)
\label{a-approximation}
\end{equation}
where 
$\phi_k$ are the usual piecewise linear hat functions defined over a partition of $[c_{\min},c_{\max}].$ Note that any function $a(c)$ can be represented by its corresponding vector ${\bf a}.$
%
 Since the values of $c_{\min}$ and $c_{\max}$ may vary considerably for each $a$ used during the optimization, we choose instead an interval that is sufficiently large to include the range of $c$ for each $a$ considered by the algorithm. In practice, this means making a guess, and expanding the interval when $c$ leaves our chosen interval.

The penalty term $\alpha\norm{a-a^*}^2_{\mathcal{A}}$ can be replaced by
$$\norm{\sum_{i=1}^L (a_i-a^*_i)\phi_i}_{\mathcal{A}} =({\bf a-a^*})^T B ({\bf a-a^*})$$
where
the components of the matrix $B$ are given by
$B_{ij}= (\phi_i,\phi_j)_{\mathcal{A}}.$

Various strategies for the choice of regularization parameters are discussed in \cite{Vogel}. In each of the following examples, we 
used an L-curve method to choose an optimal regularization parameter $\alpha.$

Recall that our chemotaxis system is
\begin{eqnarray} 
&u_t = M\Delta u - \nabla \cdot (a(c)u\nabla c) & \mbox{in $\Omega \times (0,T)$} 
\\ \nonumber 
&c_t= D\Delta c + {bu\over u+h} -  \mu c &\\ \nonumber 
&u(x,0) =u_0(x), \quad c(x,0) = c_0(x) & \mbox{for $x\in \Omega$}\\ 
&\frac{\partial u}{\partial \nu}= \frac{\partial c}{\partial \nu}= 0 & \mbox{on $\partial\Omega\times (0,T).$}
\nonumber 
\end{eqnarray}
A similar system was used by Myerscough et al. \cite{Myerscough} in their numerical simulations of chemotaxis in limb-bud development 
with parameters
$$ M=0.25, D=1,  a(c)=2, h=1, b=\mu, u_0=1+\varepsilon(x) , c_0=0.5, \Omega=[0,1]$$
where $\mu \in [0,3000]$ and $\varepsilon(x)$ was a bounded perturbation function.
In Examples 1-3, we used the Myerscough parameters with
$$\varepsilon(x) = e^{-55(x-0.5)^2}, T=0.25, b=\mu=50.$$

\begin{figure}[hp]
\begin{center}
\resizebox{5in}{!}{\includegraphics{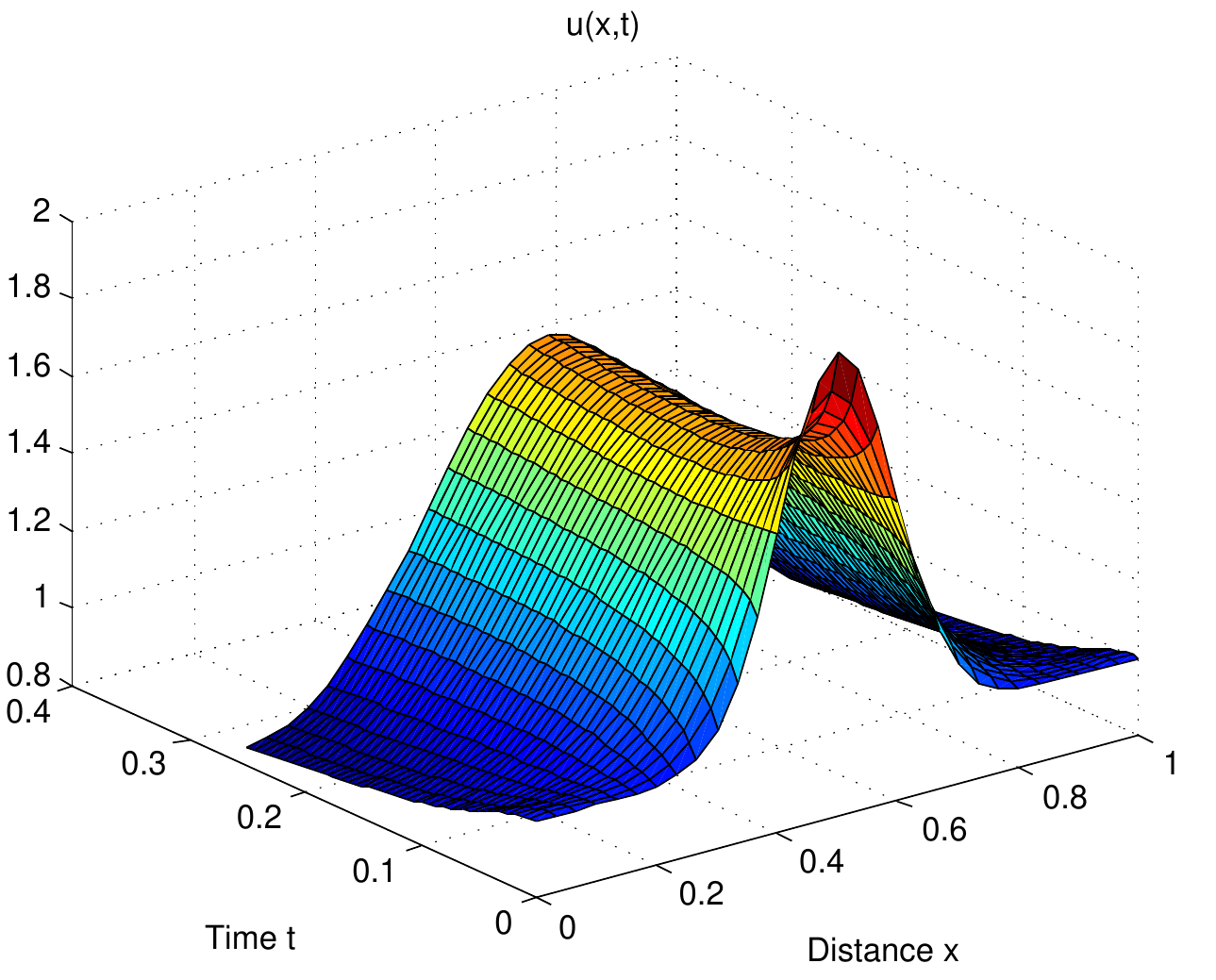}\hspace{1in}\includegraphics{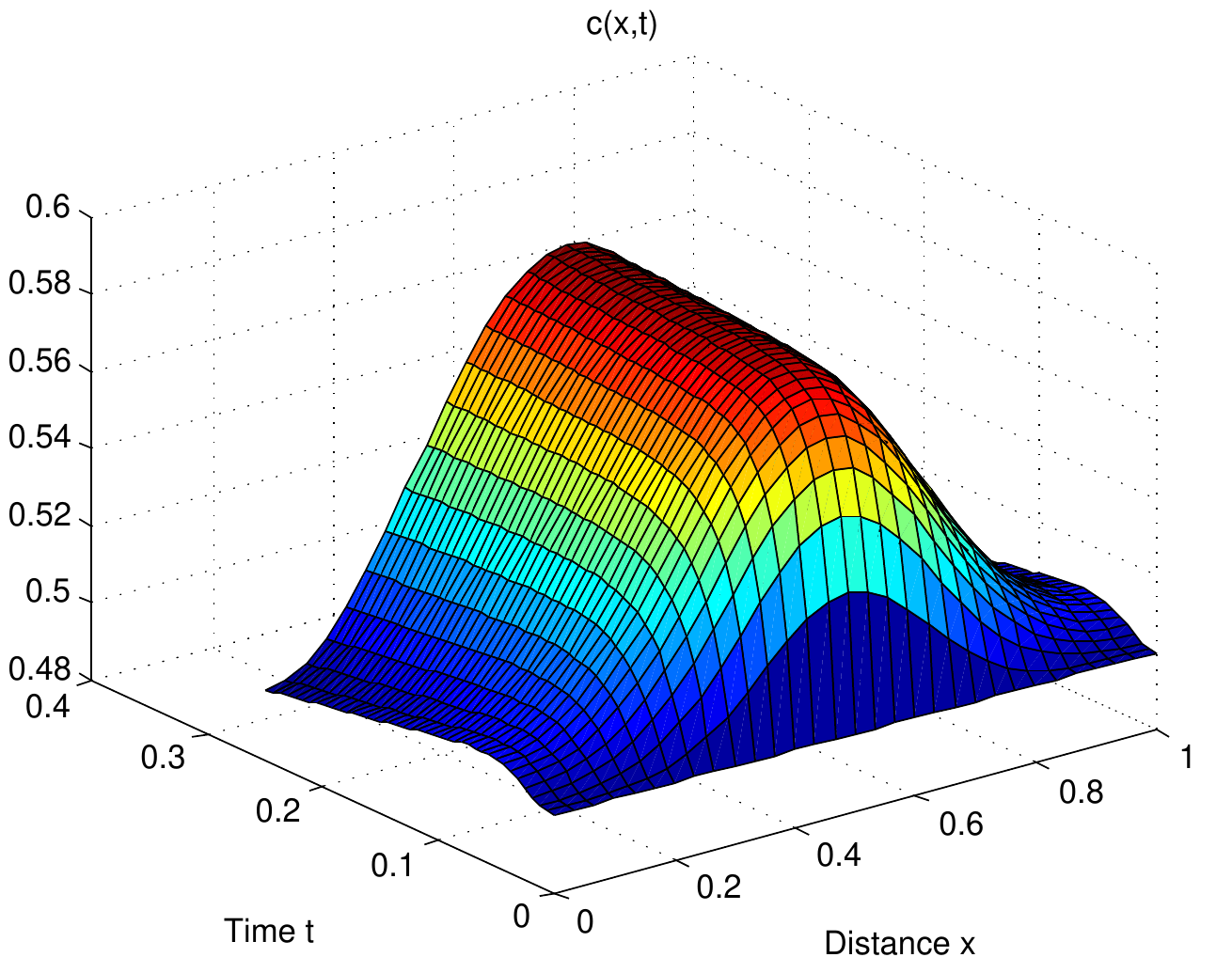}}
\end{center}
\caption{Cell and Chemical Concentrations over time interval $[0,0.25]$ with $a(c,u,x,t)=2.$}
\label{concentrations_example1}
\end{figure}

\paragraph{Example 1}
Consider the chemotactic coefficient used by Myerscough et al. \cite{Myerscough}
$$ a(c,u,x,t)= 2.$$
The cell and chemoattractant concentrations associated with this $a$ are
shown over the time interval  $[0,0.25]$
in  Figure \ref{concentrations_example1}. 
An initial guess of $a=1$ was used. The {\em a priori} guess was also chosen to be $a^*=1.$
The parameter $a$ was recovered to within 
$1.461 \times 10^{-6}$ of the true value at $T=0.25.$
%
Commented these next two examples out!

\begin{figure}[h]
\begin{center}
\resizebox{5in}{!}{\includegraphics{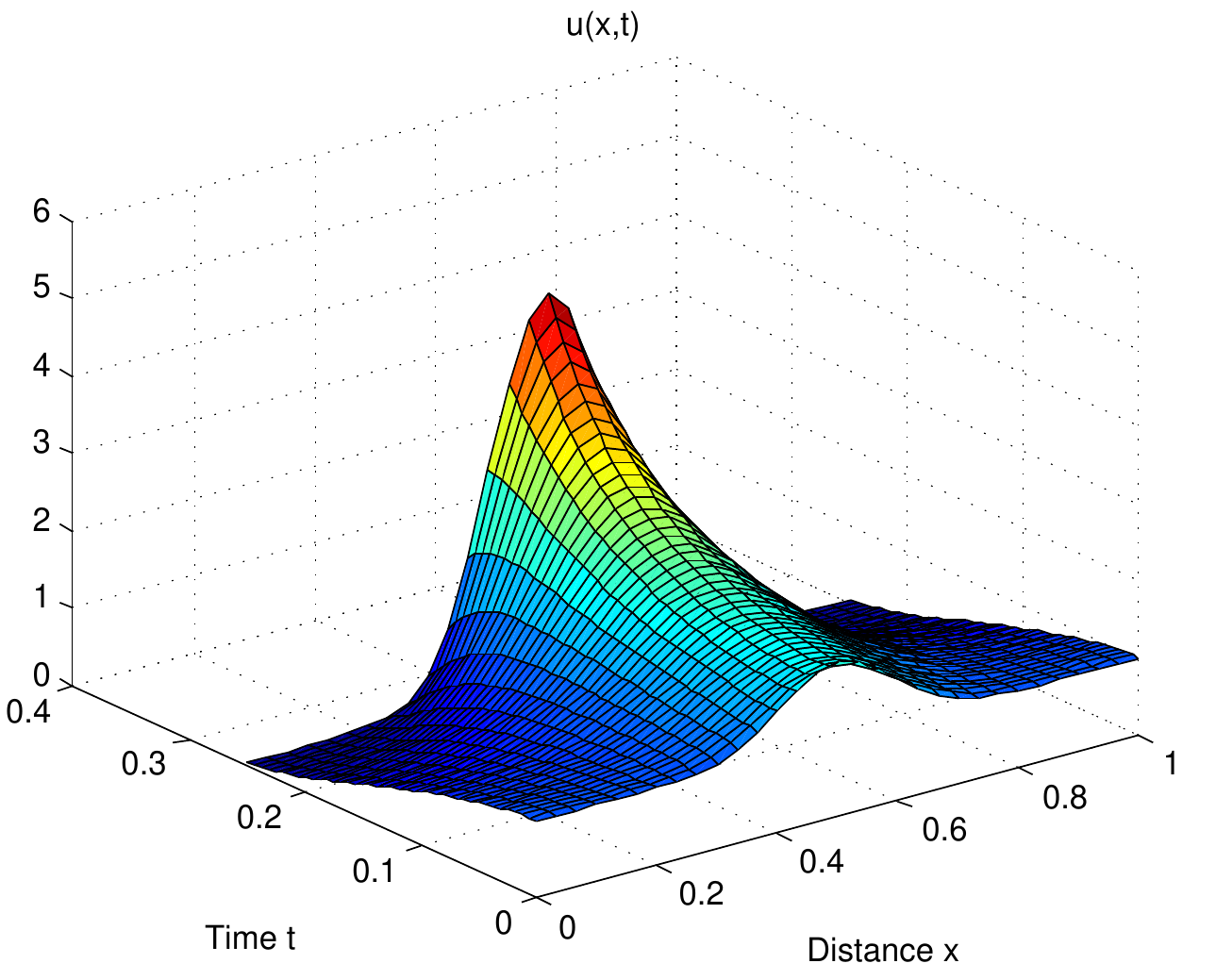}\hspace{1in}  \includegraphics{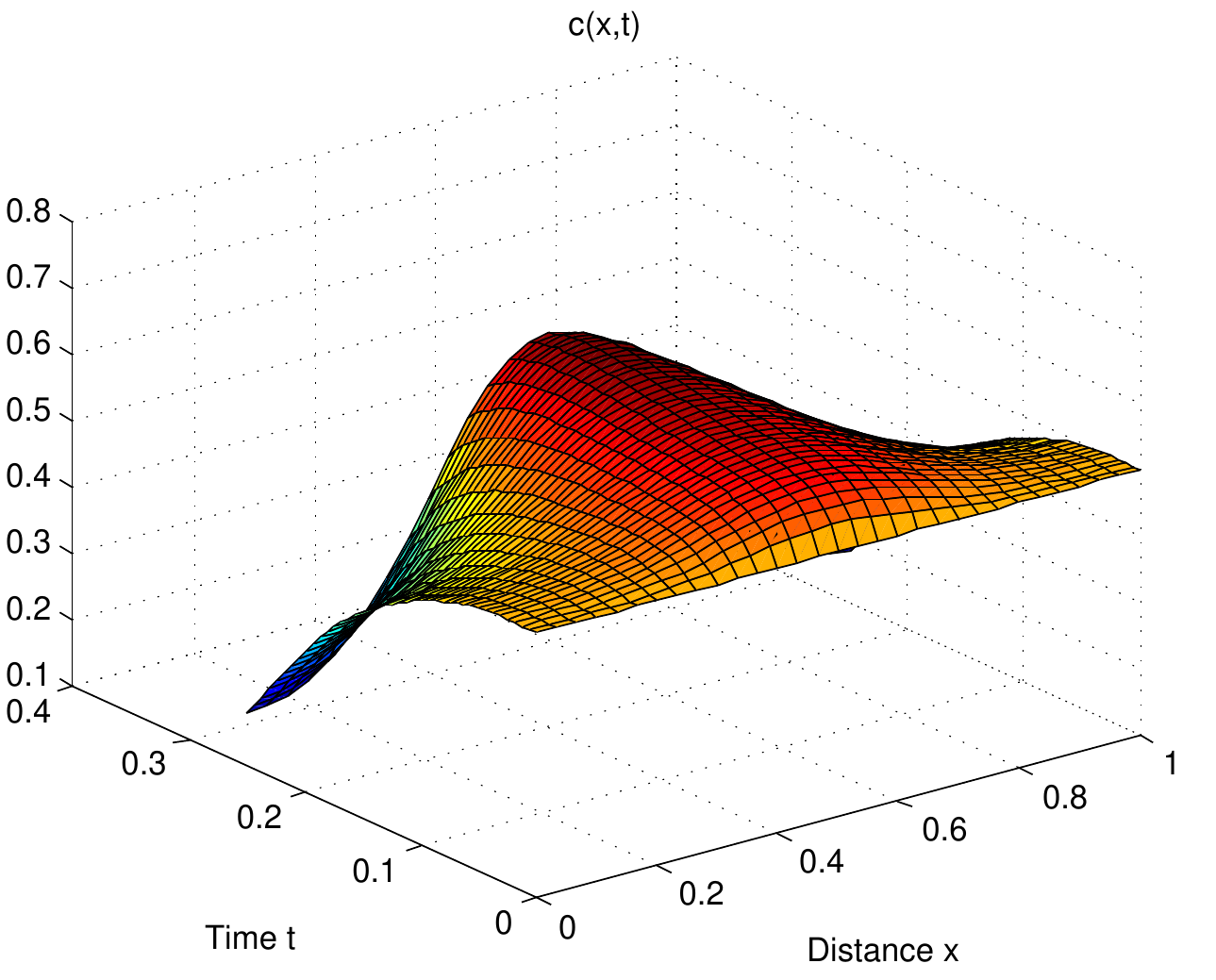}}
\end{center}
\caption{Cell and Chemical Concentrations over time interval $[0,0.25]$ with $a(c,u,x,t)= 2/c$ }
\label{concentrations_example4}
\end{figure}

\paragraph{Example 2: Keller-Segel model}
We consider the nonlinear chemotactic coefficient
$$ a(c)= 2/c$$
proposed in the original Keller-Segel model for chemotaxis \cite{Keller}.
The cell and chemoattractant concentrations associated with this $a$ are
shown at steady state in Figure \ref{concentrations_example4}. 

From the data, we find that  $[c_{\min},c_{\max}] = [0.1794, 0.6398]$  when $a(c)=2/c.$ 
Since the optimization algorithm will use approximations of other chemotactic functions, we attempt to reconstruct $a(c)$ over a larger interval.
We found the interval  $[0.1,0.7]$ to be sufficiently large to include the range of $c$ for each $a$ considered by the algorithm. 

\begin{figure}[h]
\begin{center}
\resizebox{5in}{!}{\includegraphics{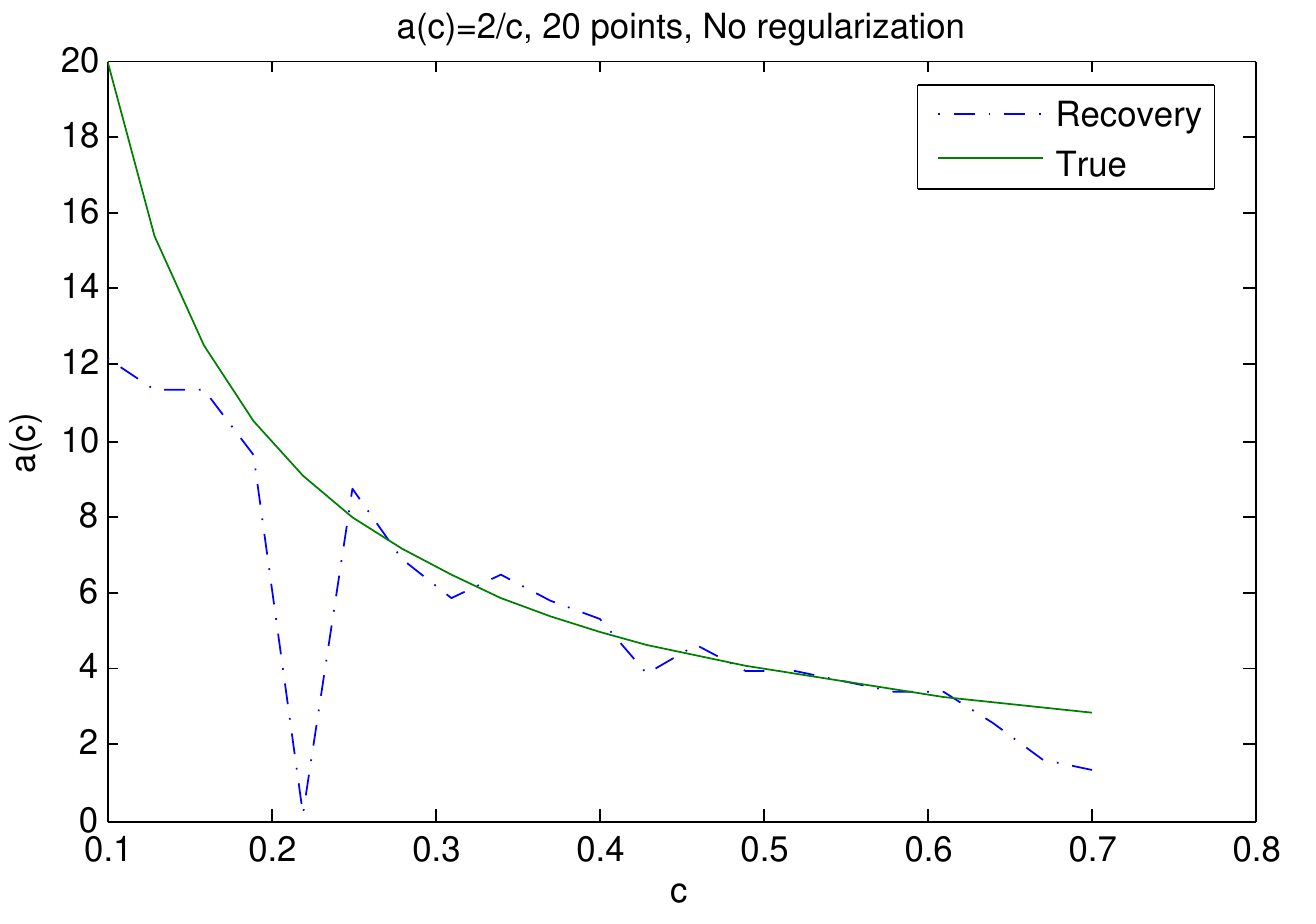} \hspace{1in}  \includegraphics{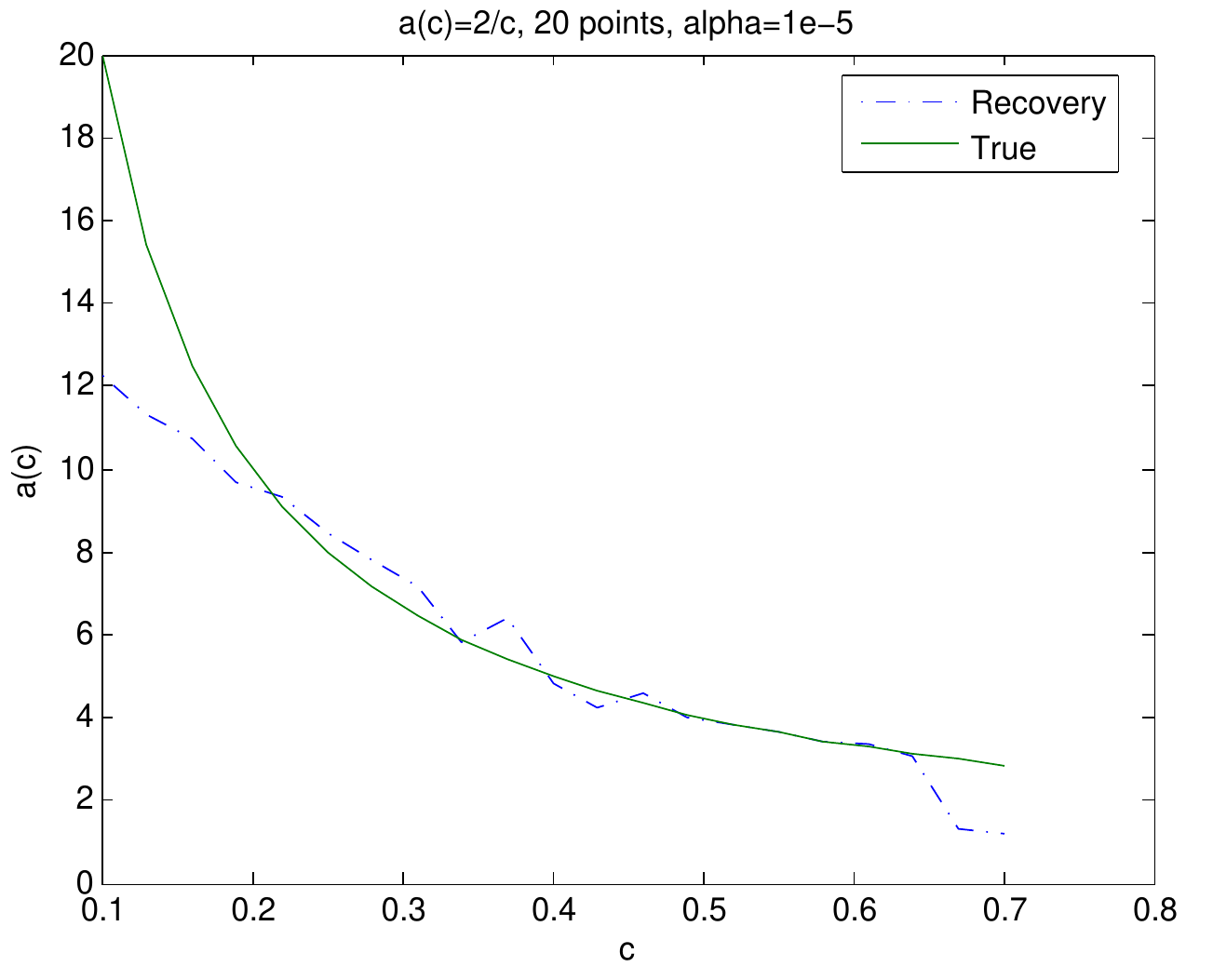}}
\end{center}
\caption{True chemotactic coefficient $a(x) =2/c$ and its recovery with no regularization (left) and with 
a regularization parameter of $\alpha=10^{-5}.$}
\label{a_example4}
\end{figure}

An initial guess of $a=15(1-c)^2$ was used.
The {\em a priori} guess was also chosen to be $a^*=15(1-c)^2.$
Figure \ref{a_example4} shows the chemotactic function $a$ and its recovery of $a_{noise}$ with and without regularization. 
A regularization parameter of  $\alpha=10^{-5}$ was chosen by an L-curve method. Notice that the regularized recovery is quite reasonable over the interval $[c_{\min},c_{\max}] = [0.1794, 0.6398]$ and that its quality degrades, as expected, outside this interval.

\paragraph{Comments}
The number of iterations used by the algorithm is quite sensitive to the choice of initial function $a_0$ and the number of piecewise linear basis functions used in  equation (\ref{a-approximation}). For experimental data, we must acknowledge that the quality of the recovery degrades with increased noise in the data. In certain applications such as pattern formation in {\em Escherichia coli} or {\em Salmonella typhimurium}, see Tyson et al. \cite{Tyson}, the size of the interval $[c_{\min},c_{\max}]$ is sometimes too small to give adequate information for the recovery of the chemotactic coefficient. This can be avoided by taking a larger time interval $[0,T].$ In numerical simulations, this requires a careful choice of numerical method for the solution of the chemotaxis system, see \cite{Tyson2}. An alternative approach is to restrict our measurements to a particular time, rather than an interval of time.  The efficacy of this approach will be discussed in a future paper.

\section{Conclusions}
In this work, we have explored a particular mathematical aspect of the chemotactic sensitivity within the gradient. The identification of a chemotactic sensitivity with functional dependence has been determined. The interesting aspect of this work is that, to our understanding, no one has been able to capture the chemotactic sensitivity information from limited data with dependence on the chemical in a system. We have proven the existence of the state solutions in specific Sobolev spaces and formulated an inverse problem. We have employed Tikhonov regularization to recover the chemotactic sensitivity from noisy measurements. In doing so, a minimization problem is formed and the necessary convergence results for an approximating minimizer to the true parameter are discussed. 

Another significant result is that we have established a source condition that guarantees a particular rate of convergence by imposing a Lipschitz condition on the derivative of the chemotactic sensitivity. In practice, this is biologically reasonable, since the chemotactic sensitivity has a rate of change that is bounded for bacterial growth, \cite{Ford}.

Numerically, we have utilized models from
Myerscough et. al, \cite{Myerscough} and Keller and Segel, \cite{Keller}  for the studies of the comparison of our proposed work to the actual scenarios. 
With the use of Tikhonov regularization, we have been able to recover the chemotactic sensitivity with reasonable accuracy. A biological benefit of this knowledge is the ability for one to understand the growth associated with chemotaxis within tumor studies, leukocyte dynamics, and bacterial patterns based on the specific gradient information that can be recovered from imperfect data.

\section{Acknowledgements}

This work was partially supported by  the Howard Hughes Medical Insitute as part of its Undergraduate Biological Sciences Education Program award to Murray State University, and by the National Science Foundation awards DMS-0209562, DMS-0414011 and DUE-0531865. Any opinions, findings, and conclusions or recommendations expressed in this material are those of the authors and do not necessarily reflect the views of the  Howard Hughes Medical Insitute or the National Science Foundation.

\appendix
\section{Appendix}
\subsection{Yagi's existence theory}
For more thorough understanding of Theorem \ref{existencet}, we include the parabolic problem convention used in Yagi, \cite{Yagi}.\begin{eqnarray} 
&u_t = \nabla \cdot (a(u,p)\nabla u)-ub(p) \nabla p) & \mbox{in $\Omega \times (0,\infty),$} 
\label{Yagiprob}
\\ \nonumber
&p_t= d\Delta p +uf(p)-g(p)p & \mbox{in $\Omega \times (0,\infty),$} 
\\ \nonumber 
&u(x,0) =u_0(x), \quad p(x,0) = p_0(x) & \mbox{for $x\in \Omega$}\\ 
&\frac{\partial u}{\partial \nu}= \frac{\partial p}{\partial \nu}= 0 & \mbox{on $\partial\Omega\times (0,\infty).$}
\nonumber 
\end{eqnarray}In Yagi's theorem \cite[Thm 2.1]{Yagi}, it states \begin{theorem} 
\label{existencelocalYagi} 
Let $u_0, p_0 \in H^{1+\ve_0} (\Omega)$ for $0<\ve_0\le 1$, 
and $u_0(x)\ge 0$, $p_0(x) \ge > \delta_0>0$ on $\overline \Omega$.  Assume a real local solution $(u,p)$ to  (\ref{Yagiprob}) exists on the interval $[0,S]$ such that 
$$ 
u, p \in C \([0,S) ; H^{1+\ve_1} (\Omega) \) \cap C \([0,S) ; H^2 
(\Omega) \) \cap C^1 \([0,S) ; L^2 (\Omega) \)$$
for some $\ve_1>0.$ In addition, assume that $p$ satisfies $p(x,t)>0$ on $\overline{\Omega} \times [0,S]$ and an estimate$$\norm {p(t)}_{H^2}\le At^{{\frac{(\ve_{2}-1)}{2}}} on \,{\mbox 0<t\le S,} $$ for some $\ve_{2}>0$ and constant $A$. Then  $$ 
u(x,t) \ge 0, \quad p(x,t) \ge \underline{p}(t)\quad \mbox{ for all $(x,t) \in \overline{\Omega}\times[0,S]$}, $$ where $\underline{p}$ denotes a positive function defined as the global solution to an ordinary differential equation,\begin{eqnarray}\label{peqn}\frac {d{\underline{p}}}{dt}&=&-g(\underline{p})\underline{p}\,  {\mbox { on $0<t< \infty$}}\\ \nonumber\underline{p}(0)&=&\delta_0>0.\\ \nonumber\end{eqnarray}

\end{theorem}

For further connection to our work, we utlilize the continuation to a unique solution that Yagi developed  \cite[Thm 3.4]{Yagi} with\begin{theorem}\label{uniqueYagi}Let $u_0, p_0 \in H^{1+\ve_0} (\Omega)$ for $0<\ve_0\le 1$, 
and $u_0(x)\ge 0$, $p_0(x) \ge \delta_0>0$ on $\overline \Omega$ and let $0<\eta<\beta-\alpha$.  Then, in the function space, $C ^{\eta}\([0,\infty) ; H^{1+\ve_1} (\Omega) \)$, the problem (\ref{Yagiprob}) possess a unique local solution$$u, p \in C \([0,S) ; H^2 
(\Omega) \) \cap C^1 \([0,S) ; L^2 (\Omega) \)$$with the lower bounds$$ 
u(t) \ge 0, \quad p(t) \ge \underline{p}(t)\quad \mbox{ for $t\in[0,S]$}, $$ where $\underline{p}(\cdot)$ is a positive function defined by equations in (\ref{peqn}).  The interval $[0,S]$ on which the solution exists at least is determined by the norms $\norm{u_0}_{H^{1+\ve_0}}$, $\norm{p_0}_{H^{1+\ve_0}}$ and by the initial lower bound $\delta_0$.\end{theorem}
It is to be noted that by this theorem from Yagi's work 
that a maximal solution to (\ref{Yagiprob}) can be uniquely defined in the space $C^{\eta} \([0,S) ; H^{1+\ve_1} (\Omega) \)$ for $0<\eta<\beta-\alpha$ for each $u_0, p_0$ such that  $u_0, p_0 \in H^{1+\ve_0} (\Omega)$ for $0<\ve_0\le 1$, 
and $u_0(x)\ge 0$, $p_0(x) \ge \overline c_0 > \delta_0>0$ on $\overline \Omega$.

\end{document}